\documentclass[11pt]{amsart}
\usepackage{amsfonts,amsmath,amssymb,comment,amscd,latexsym,}
\numberwithin{equation}{section}

\newtheorem{Lem}{Lemma}[section]
\newtheorem{Prop}[Lem]{Proposition}
\newtheorem{Cor}[Lem]{Corollary}
\newtheorem{Thm}[Lem]{Theorem}

\theoremstyle{definition}
\newtheorem{Def}[Lem]{Definition}

\renewcommand{\k}{{\Bbbk}}

\newcommand\e{e_{\lambda}}
\newcommand\ex{e_{{\lambda}X}}
\newcommand\R{R_{\kappa,N}}

\newcommand{\G}{{\widehat G}}
\newcommand{\D}{\mathcal{D}}
\newcommand{\pf}{\medskip\noindent{\sc Proof:\,}}

\begin{document}

\title[Representation theory of liftings]{Representation theory of liftings of quantum planes}
\author{William Chin}
\address{DePaul University\\Department of Mathematical Sciences\\2320 N. Kenmore\\Chicago, IL 60614}\email{wchin@condor.depaul.edu} 
\author{Leonid Krop}
\address{DePaul University\\Department of Mathematical Sciences\\2320 N. Kenmore\\Chicago, IL 60614}\email{lkrop@depaul.edu}

\begin{abstract}We determine the regular representations, Gabriel quivers and representation type of all liftings of two-dimensional quantum linear spaces.
\end{abstract}
\date{10/17/08}
\maketitle

\section*{Introduction}\label{intro}

Liftings of quantum linear spaces were constructed and completely described in \cite{AS1}, and independently in \cite{BDG}. These Hopf algebras belong to the class of pointed finite-dimensional Hopf algebras with an abelian group of group-likes and are, arguably, the simplest Hopf algebras constructed by the lifting method of N. Andruskiewitsch and H.-J. Schneider (see e.g. \cite {AS2}). When the dimension of a quantum linear space is two, we colloquially refer to them as quantum planes. In this article we systematically study the representation theory of liftings of quantum planes. The group of group-likes is arbitrary finite abelian in our construction, and the base field $\Bbbk$ is algebraically closed of characteristic zero. Liftings of quantum planes are generated as algebras by two skew-primitive elements along with the group-likes $G$. The skew-primitives generalize the generators $e,f$ in the restricted enveloping algebra of $\mathfrak sl(2)$ and quantum analogs at roots of unity. If there is a nontrivial commutation relation between the skew-primitives we say that the lifting is linked.

Let us summarize some existing related work. The simple representations of liftings of quantum planes were treated in \cite{{AB1},{AB2}} where simple modules are described or reduced to known theory in most cases. The representation theory of various versions of quantized restricted enveloping algebras was studied in \cite{{CP}, {XJ}, {SU}} using a variety of techniques. In \cite{EGST} representations of quantum doubles of generalized Taft algebras were examined. The simple and projective modules were explicitly constructed within the regular representation, extending methods from \cite{SU}. These Hopf algebras are examples of liftings of quantum planes in the linked nilpotent case. In \cite{KR} rank one Hopf algebras were constructed and the structure of the regular representation of their doubles was obtained. These Hopf algebras are sometimes liftings of quantum planes; but this is not always the case as sometimes they are not even pointed. 

In this article we explicitly describe the simples, projectives, blocks, and Gabriel quiver in all cases. This enables us to determine the representation types of the blocks. We make use of certain central idempotents that are constructed from certain equivalence classes of characters of $G$ to reduce to the corresponding two-sided summands, which we call {\em class subalgebras}. These idempotents were also used in \cite{{AB1},{AB2}}. The class subalgebras we encounter are generated by three elements, one of which is a unit, and the other two elements are the images of skew-primitives. Our analysis breaks into the consideration of summands where both, just one, or neither of the generators is nilpotent. We shall refer to these cases as {\em nilpotent, seminilpotent} and {\em unipotent}, respectively. Another division is into the classes of linked and unlinked liftings, resulting in six distinct cases altogether.

Let us give an outline of the results. In the unlinked nonnilpotent cases, the class subalgebras are the blocks, and are either Nakayama algebras, or skew group algebras over truncated polynomial rings. The Nakayama algebras are the simplest algebras of finite representation type, while the truncated polynomial rings are among simplest of the algebras of wild and tame representation type (depending on the degree of truncation). In the linked nonnilpotent cases the class subalgebras are direct sums matrix algebras over a Nakayama algebras and are thus of finite representation type.

In the remaining linked nilpotent case the nonsimple blocks are special biserial algebras and therefore of tame representation type. The algebras in this case generalize the quantized restricted enveloping algebra of $\mathfrak sl(2)$. The basic algebras and quivers that occur here also arose in the study of Hopf algebras in \cite{EGST}, and we can apply their results to this case. Here we use an alternative less explicit but much simpler determination of the simple and projective representations using analogs of ``baby Verma modules'' in the more general setting of seminilpotent and nilpotent liftings of quantum planes. We show that linked liftings are symmetric algebras and we provide an analog of the Casimir element. These tools allow a concise and efficient determination of the socle and Loewy series structure of the projective indecomposables, and the minimal central idempotents of nonsimple blocks.

A more detailed description of material by sections is as follows. In section \ref{prelims} we review the ideas leading up to the construction of liftings of quantum linear spaces. We define types of liftings and discuss some auxiliary facts about duality for finite abelian groups and group algebras used in the sequel.

In section \ref{Gen} we introduce the general theory of liftings of quantum planes. We define class idempotents, class subalgebras and compute bases for them. A quick argument proving Lemma \ref{centralid} shows that the class idempotents are central. In \S 2.1 we give generators and relations for liftings and comment on the existence of unlinked or linked data.

In section \ref{unlinked} we determine the structure of simples, projectives, Gabriel quiver and blocks for unlinked liftings. The nilpotent case is handled in Theorem \ref{nilunlink}, where the class subalgebras are the blocks and are essentially skew group rings over truncated polynomial rings. The seminilpotent case is described in Theorem \ref{seminilpunlinked} where the class subalgebras either resemble the nilpotent case or are nonbasic Nakayama algebras. In this case the Gabriel quiver has vertices corresponding to certain cosets of a certain subgroup of the character group. In the unipotent case, we show in Theorem \ref{unipotunlink} that the class subalgebras are semisimple, and we explicitly construct the simple modules. We also give a criterion for the isomorphism of the simples in Theorem \ref{unipotent}.

We turn to the linked liftings in section \ref{linked}, where we find an analog of the Casimir element from classical Lie theory, which is central by Lemma \ref{cas}. In Proposition \ref{symalg} we show that linked liftings are symmetric algebras by directly showing unimodularity and that the square of the antipode is inner. In \S 4.2 we study the image of the Casimir element C in each class subalgebra and compute the minimal polynomials. In the unipotent case we quickly find in Thereom \ref{potentclsalg} that each class subalgebra is isomorphic to a matrix algebra over the subalgebra generated by C. As a result, each class subalgebra is a direct sum of matrix rings over local rings. In the seminilpotent case, addressed in Theorem \ref{lsn} where the blocks are matrix rings over the base field or over the truncated polynomial algebra $\Bbbk[v]/(v^2)$. The precise decompositions are expressed in a number of cases.  For use in the  seminilpotent and  nilpotent  cases,  we  introduce analogs of standard cyclic modules in \S 4.5, where they are shown to be simple in Proposition \ref{simpleinduced} in the seminilpotent case. Finally we study the lengthier nilpotent case in \S 4.6.  Here the standard cyclic modules may no longer be simple, but have a unique maximal submodule. Paralleling Lie theory, we obtain precise results concerning the simple quotients. As a first step in Theorem \ref{dumbcase}, we dispose of the generic case where the class subalgebra is semisimple. We then look at the complementary case and we give detailed results about the Loewy factors of the standard cyclic modules. This enables the determination of the structure of the projective indecomposable modules in Theorem \ref{main}, along with a presentation of the basic algebras of blocks by quivers with relations in Theorem \ref{quiver}. We close by explicitly giving embeddings of the projective indecomposables into nonsimple blocks.

\section{Preliminaries}\label{prelims}

 {\bf Notation}: $G$ a finite abelian group

$\G$ the character group of $G$

$\Bbbk$ an algebraically closed field

$a_1,\ldots,a_n\in G$

$\chi_1,\ldots,\chi_n\in\G$

$q_{ij}=\chi_j(a_i),\quad q_i=\chi_i(a_i)$

For $g\in G$ $|g|$ denotes the order of $g$. In particular, for $q\in\Bbbk^{\bullet}$ $|q|$ is the order of $q$.

For $\kappa\in\Bbbk$ and $N\in\Bbb N$,\quad $R_{\kappa,N}=\{\gamma\in\Bbbk|\gamma^N=\kappa\}$
\subsection{Liftings of Quantum Linear Spaces}\label{liftings}

Recall \cite{{Yet},{RT},{AS2}} that the category ${}^G_G\mathcal {YD}$ consists of $\Bbbk G$-modules and $\Bbbk G$-comodules $V$ such that the 
$G$-grading 
$$V=\bigoplus V_g\;\text{where}\;V_g=\{v\in V|\rho(v)=g\otimes v\}$$
where $\rho:V\to\Bbbk G\otimes V$ is the comodule structure map,  satisfies 
$$h.V_g\subset V_{hgh^{-1}}\;\text{ for every}\; h,g\in G.$$
Since $G$ is abelian $h.V_g=V_g$ and $V$ has a basis, say, $\{v_i|i\in I\}$ of $G$ and $\G$- eigenvectors indexed by a set $I$. In other words,
\begin{align}g.v_i&=\chi_i(g)v_i\label{bihom1}\\ 
                                 \rho(v_i)&=a_i\otimes v_i\label{bihom2}\end{align}
for some $a_i\in G$ and $\chi_i\in\G$, and for all $g\in G,\;i\in I$.

We say that an element $v\in V$ is {\em bihomogeneous} of degree $(a,\chi)$ if the equations \eqref{bihom1} and \eqref{bihom2} hold for $v,g$ and $\chi$ and we write $a=g_v,\chi=\chi_v$. 
\begin{Def}(\cite{AS1}) An $n$- dimensional Yetter- Drinfel'd module $V$ is called a quantum linear space if
$$\chi_i(a_j)\chi_j(a_i)=1\;\text{for all}\;i\neq j$$
in some bihomogeneous basis for $V$.
\end{Def}

We review a construction of the Nichols algebra $B(V)$ associated to the quantum linear space $V$ (cf. \cite[Lemma 3.4]{AS1}). As an algebra $B(V)$ is defined via the relations
\begin{align*}v_iv_j&=\chi_j(a_i)v_jv_i\\
v^{n_i}&=0\end{align*}
for all $i\neq j$ and $n_i=|q_i|$. One can see immediately that the set $\{v_1^{i_1}\cdots v_n^{i_n}\}$ is a basis of $B(V)$. The action and coaction of $G$ on $V$ both extend uniquely to $B(V)$ by requiring $B(V)$ to be an algebra in ${}^G_G\mathcal{YD}$. Explicitly
\begin{align*}g.v_1^{i_1}\cdots v_n^{i_n}&=\chi_1^{i_1}\cdots\chi_n^{i_n}(g)v_1^{i_1}\cdots v_n^{i_n}\\
\rho(v_1^{i_1}\cdots v_n^{i_n})&=a_1^{i_1}\cdots a_n^{i_n}\otimes v_1^{i_1}\cdots v_n^{i_n}.\end{align*}
Note that every monomial $v_1^{i_1}\cdots v_n^{i_n}$ is bihomogeneous.

The coalgebra structure  maps $\delta,\epsilon$ are defined by first setting
$$\epsilon(v_1^{i_1}\cdots v_n^{i_n})=\delta_{0,i_1}\cdots\delta_{0,i_n}$$
and for $\delta:B(V)\to B(V)\otimes B(V)$ is given by
$$\delta(v_i)=v_i\otimes 1+1\otimes v_i$$
on generators and extends to $B(V)$ as follows. Following Lusztig \cite{L} or \cite{AS2} we give $B(V)\otimes B(V)$ a new multiplication by setting
$$(r\otimes s)(t\otimes u)=\chi_t(g_s)rt\otimes su$$
for monomials $r,s,t,u$, and denote this algebra by $B(V)\underline{\otimes} B(V)$. We also endow $B(V)\otimes B(V)$ with the usual tensor product module and comodule structures over $\Bbbk G$. A tedious, but straightforward calculation shows that $B(V)\underline{\otimes} B(V)$ is an algebra in ${}^G_G\mathcal{YD}$, $\delta$ is well-defined and $\delta,\epsilon$ are $G$- linear and $G$- colinear algebra maps. The above mentioned properties amount to saying that $B(V)$ is a bialgebra in ${}^G_G\mathcal {YD}$. We turn $B(V)$ into a Hopf algebra in ${}^G_G\mathcal {YD}$ by defining an antipode $\mathcal S$ as the linear map
$${\mathcal S}(v_1^{i_1}\cdots v_n^{i_n})=(-1)^{\sum i_j}\prod_{j=1}^nq_j^{\binom{i_j}2}v_1^{i_1}\cdots v_n^{i_n}.$$ 
We can use now a result of Radford \cite{RA}, recast by Majid \cite{MA} in categorical terms, to the effect that the biproduct or bosonization $B(V)\#\Bbbk G$ is a Hopf algebra. The algebra $B(V)\#\Bbbk G$ is is known as the {\em trivial lifting} of $V$.
\begin{Def} A Hopf algebra $H$ is a lifting of $V$ if the graded Hopf algebra associated to the coradical filtration of $H$ is isomorphic to $B(V)\#\Bbbk G$.\end{Def}

Next we derive a general property of liftings.
\begin{Prop}Every lifting of $V$ has a structure of left-left Yetter-Drinfel'd module.
\end{Prop}
\pf Let $H$ be a lifting of $V$. By \cite{AS1} $H_0=\Bbbk G$ hence $H$ is naturally a $G$-module under the action of $G$ by conjugation,
 $h\mapsto g.h:=ghg^{-1}$ for every $h\in H,g\in G$. Moreover, $H$ is generated by $G$ and skew-primitives $x_i$ satisfying $g.x_i=\chi_i(g)x_i$. By \cite[5.2]{AS1} the set $\{gx_1^{i_1}\cdots x_n^{i_n}|g\in G,0\le i_j<n_j\}$ is a basis for $H$. Let $I$ be the span of all $gx_1^{i_1}\cdots x_n^{i_n}$ with $\sum i_j>0$. By the quantum binomial formula \cite{KA} one can see readily that $I$ is a coideal and it complements $H_0$. Therefore we have a coalgebra projection $\pi:H\to\Bbbk G$. Then $\rho:=(\pi\otimes\text{id})\Delta: H\to \Bbbk G\otimes H$ equips $H$ with a $\Bbbk G$- comodule structure. Further, every $u=gx_1^{i_1}\cdots x_n^{i_n}$ is bihomogeneous of degree $(ga_1^{i_1}\cdots a_n^{i_n}, \chi_1^{i_1}\cdots \chi_n^{i_n})$. Therefore for every $h\in G$
$$\rho(h.u)= \chi_1^{i_1}\cdots \chi_n^{i_n}(h)\rho(u)=ga_1^{i_1}\cdots a_n^{i_n}\otimes h.u$$
and the proof is complete.\qed

We apply the previous proposition and consider the braided commutator $[\,,\,]_c:H\otimes H\to H$ corresponding to the braiding $c:H\otimes H\to H\otimes H$ arising from the Yetter- Drinfel'd module structure on $H$ \cite{AS2}. This commutator is given by $[a,b]_c=\mu(\text{id}-c)(a\otimes b)$ for $a,b\in H$ where $\mu$ is the multiplication in $H$. When $a,b$ are bihomogeneous we have
$$[a,b]_c=ab-\chi_b(g_a)ba. $$

We shall later have use of the following braided commutator rules. Let $x,y$ be bihomogeneous elements of degrees $(a,\chi)$ and $(b,\chi^{-1})$, respectively, satisfying
\begin{equation}\label{simplecommutator}[x,y]_c=ab-1\quad\text{and}\quad\chi(a)=\chi(b).\end{equation}
\begin{Prop}\label{productrule}Let $H$ be a lifting of $V$. Suppose $u,v\in H$ are bihomogeneous and $x,y$ are as above. Then
\begin{itemize}
\item[(a)] $[x,uv]_c=[x,u]_cv+\chi_u(a)u[x,v]_c$
\item[(b)] Let $q=\chi^{-1}(a)$. For every $s\ge 1$
\begin{itemize}
\item[(i)] $[x,y^s]_c=(s)_qy^{s-1}(q^{s-1}ab-1)$
\item[(ii)] $[y,x^s]_c=q^{-s}(s)_qx^{-s+1}(q^{-(s-1)}ab-1)$.
\end{itemize}
\end{itemize}
\end{Prop}
\pf (a) is seen by a direct inspection.

(b) (i). The formula holds for $s=1$ by \ref{simplecommutator}. We induct on $s$ assuming it holds for the given $s$. By part (a) we calculate
$$[x,y^{s+1}]_c=[x,y^s\cdot y]_c=(s)_qy^{s-1}(q^{s-1}ab-1)y+q^sy^s(ab-1)$$
because $y^s$ is a $G$- eigenvector of weight $\chi^{-s}$. Since $aby=q^2yab$ the  first term on the right is $(s)_qy^s(q^{s+1}ab-1)$. 
Therefore the right hand side equals
$$\gamma y^s((q^s+(s)_qq^{s+1})ab-(q^s+(s)_q))=\gamma y^s(s+1)_q(q^sab-1)$$
the last equality holds by the identities $1+(s)_qq=(s+1)_q$ and $q^s+(s)_q=(s+1)_q$.

Part (b)(ii) is proven similarly using the relation $yx-q^{-1}xy=-q^{-1}(ab-1)$ and $(s)_{q^{-1}}=q^{-s+1}(s)_q$.\qed

\subsection{Duality for abelian groups and group algebras}\label{abgroups}

We collect some basic facts that will be used freely henceforth.

Let $G$ be a finite abelian group and its dual group $\G$ as before. Then there is an inclusion-reversing correspondence between the subgroups of $G$ and the subgroups of $\G$. The correspondence takes a subgroup $H\subset G$ to
$$H^{\perp}=\{\lambda\in\G|\lambda(h)=1\,\text{for all}\,h\in H\}.$$
The inverse is defined similarly using the identification $\widehat{\G}=G$ and we note that $H^{\perp\perp}=H$. We have $(H\cap K)^{\perp}=H^{\perp}+K^{\perp}$ and $(H+K)^{\perp}=H^{\perp}\cap K^{\perp}$ \cite{MH}.

Now considering the Hopf group-algebra $\k G$ we extend the natural pairing
$$G\times\G\to\k^{\bullet},\;\langle g,\lambda\rangle\to \lambda(g)$$
to a bilinear pairing
$$\k G\times\k\G\to\k,\;\langle\sum r_ig_i,\sum\rho_i\lambda_i\rangle\to \sum r_i\rho_j\lambda_j(g_i).$$
To every $\lambda\in\G$ we associate a minimal idempotent
$$\e=\frac{1}{|G|}\sum_{g\in G}\lambda(g^{-1})g$$
of $\k G$. The subspaces $\k\e$ afford a one-dimensional representation of $G$ with character $\lambda$, i.e. $g\e=\lambda(g)\e$. Therefore 
\begin{align*}\e e_{\mu}&=\frac{1}{|G|}(\sum_{g\in G}\lambda(g^{-1})\mu(g))e_{\mu}\\
                        &=\mu(\e)e_{\mu}.\end{align*}
By orthogonality of idempotents $\e$ we have
\begin{equation}\label{ort}\langle\e,\mu\rangle=\mu(\e)=\delta_{\lambda,\mu}.\end{equation}                        
Equality (\ref{ort}) can be interpreted as saying that basis $\{\e|\lambda\in\G\}$ of $\k G$ is dual to the standard basis $\{\mu|\mu\in\G\}$ of $\k\G$. That is to say $\e$ maps to the characteristic function $p_{\lambda}$ under the Hopf algebra isomorphism $\k G\cong (\k\G)^*$. Therefore
$$\Delta(\e)=\sum_{\lambda=\mu\nu}e_{\mu}\otimes e_{\nu}.$$

%%%%%%%%%%%%%%%%%%%%%%%%%%%%%%%%%%%%%%%%%%%%%%%%%%%%%%%%%%%%%%%%%%%%%%%%%%%%%%%%%%%%%%%%%%%%%%%%%%%%%%%%%%%%%%%%%%%%%%%%%%%%%%%%%%%%%%%%%%%%%%%%%%%%%%%
\section{General Theory of Liftings}\label{Gen}
\subsection{Types of Liftings}\label{Types}
The general lifting $H$ of $V$ is given by a {\em lifting} datum
$$\D=\{G,a,b,\epsilon_1,\epsilon_2,\chi_1,\chi_2,\gamma\}$$
where $a,b\in G,\chi_i\in\G$ and $\epsilon_i\in\{0,1\},\gamma\in\k$. $H$ is generated by $G$ and $x,y$ subject to the relations of $G$ and the following
\begin{align}x^{n_1}&=\epsilon_1(a^{n_1}-1)\\
              x^{n_2}&=\epsilon_2(b^{n_2}-1)\\
              gx&=\chi_1(g)xg\\
              gy&=\chi_2(g)yg\\
              xy &-\chi_2(a)yx=\gamma(ab-1)\end{align}
           
 The coalgebra structure of $H$ is given by
 \begin{align}\Delta(x)&=a\otimes x+x\otimes 1\\
                \Delta(y)&=b\otimes y+y\otimes 1\\
                \Delta(g)&=g\otimes g\end{align}
for all $g\in G$. The datum satisfies the conditions
\begin{align}1<n_1&=|\chi_1(a)|\\
             1<n_2&=|\chi_2(b)|\\
             \chi_1(b)\chi_2(a)&=1\\
             \chi_i^{n_i}&=\epsilon\;\text{if}\;\epsilon_i=1\\
             \chi_1\cdot\chi_2&=\epsilon\;\text{and}\; ab\neq 1\;\text{if}\;\gamma\not= 0\end{align}
             
We adopt the following terminology 
\begin{itemize}
\item A lifting is {\em linked} if $\gamma\not= 0$
\item A lifting is {\em nilpotent} if $\epsilon_i=0$ for all $i$
\item A lifting is {\em seminipotent} if $\epsilon_i=0$ for exactly one $i$
\item A lifting is {\em unipotent} if $\epsilon_i=1$ for all $i$
\end{itemize}

We comment on existence of data. Assume $G$ is a product of several cyclic groups, viz. $G=\langle g_1\rangle\times\langle g_2\rangle\times\cdots$. Suppose $|g_i|=m_i, i=1,2$ and let $\theta_i$ be a primitive root of 1 of order $m_i$. Define $\chi_i\in\G, i=1,2$ by $\chi_i(g_j)=\delta_{ij}\theta_i$. The tuple $\{G,g_1,g_2,\chi_1,\chi_2\,, 0\}$ is a lifting datum. 
If $G$ is a cyclic $p$-group generated by $g$, pick an integer $s$ such that $s^2\not\equiv 0\mod{|g|}$. Set $a=g, b=g^s$, and $\chi_1(a)=\theta, \chi_2(a)=\theta^{-s}$. Then $\{G,a,b,\chi_1,\chi_2, 0\}$ is a lifting datum.

Following \cite{AS2} we call a datum {\em linkable} if $\chi_1\chi_2=\epsilon$ and $ab\neq 1$. A linkable datum can be constructed as follows. Let $G\neq\Bbb{Z}_2$. Pick a subgroup $L$ of $G$ such that $G/L$ is cyclic of order $N$, generated by $\overline{g}=gL$ with $g\in G$. Define a character $\phi:G/L\to\Bbbk^{\bullet}$ by sending $\overline{g}$ to $\theta$, where $\theta$ is a root of $1$. Let $\chi$ be the pull-back of $\phi$ to $G$. We claim that there exists $a,b\in G$ such that $b=al, l\in G$ with $ab\neq 1$. If so, the tuple $\{G,a,b,\chi,\chi^{-1},\gamma\}$ is a linkable datum for every $\gamma\in\Bbbk^{\bullet}$. 

It remains to justify the claim about the elements $a,b$. If $L=1$, then $G$ is a cyclic group of order $>2$. Hence $g^2\neq 1$. Thus $a=g=b$ will do.

Suppose $L\neq 1$. If $g^2\notin L$, then $a=g,b=gl,l\in L$ will do. In case $g^2\in L$, say, $g^2=l_0$, then either $l_0\neq 1$, and we can set $a=g=b$, or $g^2=1$ and then $a=g,b=gl, l\neq 1$ will do.
%%%%%%%%%%%%%%%%%%%%%%%%%%%%%%%%%%%%%%%%%%%%%%%%%%%%%%%%%%%%%%%%%%%%%%%%%%%%%%%%%%%%%%%%%%%%%%%%%%%%%%%%%%%%%%%%%%%%%%%%%%%%%%%%%%%%%%%%%%%%%%%%%%%%%%%%                                                     
\subsection{Class idempotents}\label{ClassIdem}
\begin{Def}\label{CI} For each $\lambda\in\G$ we associate the idempotent $\e$ as above and let $X=\langle\chi_1,\chi_2\rangle$. Let 
$$\ex=\sum_{\mu\in\lambda X}e_{\mu}$$
for each coset $\lambda X$. 
\end{Def}
For $h\in H\ex$ we shall write $\overline h$ for its image $h\ex\in H\ex$. Let $T$ be a transversal for $X^{\perp}$ in $G$. For each subgroup $L\subset G$, we shall write $\lambda_L$ for the restriction of $\lambda$ to $L$ and put 
$$e_{{\lambda}_L}=|L|^{-1}\sum_{g\in L}\lambda(g^{-1})g\in\Bbbk L.$$

\begin{Lem}\label{centralid} $\ex$ is a central idempotent in $H$ for all $\lambda\in\G$.
\end{Lem}
\pf Observe that
\begin{align*}x\e&=|G|^{-1}\sum_{g\in G}\lambda(g^{-1})xg\\&=(|G|^{-1}\sum_{g\in G}\lambda(g^{-1})\chi_1(g^{-1})g)x\\&=e_{{\lambda}\chi_1}x\end{align*}
and similarly $y\e=e_{{\lambda}\chi_2}y$. This yields the assertion.\qed
\begin{Prop}\label{basis}The following sets 
\begin{itemize}
\item[(a)] $\{\overline{t}\overline{x}^{j_1}\overline{y}^{j_2}|0\le j_i<n_i,\;t\in T\}$ and
\item[(b)] $\{e_{\lambda\chi}\overline{x}^{j_1}\overline{y}^{j_2}|0\le j_i<n_i,\;\chi\in X\}$
\end{itemize}
are  bases of $H\ex$.
\end{Prop}
\pf We have $\e g=\lambda(g)\e$ for all $\lambda\in X$ and $g\in G$. This yields $\ex g=\lambda(g)\ex$ if $g\in X^{\perp}$. It follows directly that the set in (a) in the statement spans $H\ex$. The cardinality of the set in (a) is $[G:X^{\perp}]n_1n_2$, and since $[G:X^{\perp}]=|X|$, the set has cardinality $|X|n_1n_2$. Summing over a transversal for $X$ in $\G$ gives a spanning set for $\oplus H_{\lambda X}$ having cardinality $[\G: X]|X|n_1n_2=|G|n_1n_2$. This is just the dimension of $H$, so our spanning set is a basis of $H$. This trivially implies that the spanning set in (a) is a basis for $H\ex$.

The set in (b) is clearly a spanning set for $H\ex$, and as in the proof of (a), it is a basis because it has the requisite cardinality $|X|n_1n_2$.\qed
\begin{Lem}\label{idempotents}$\ex=e_{{\lambda}_{X^{\perp}}}$
\end{Lem}
\pf Since $e_{{\lambda}_{X^{\perp}}}\in\Bbbk X^{\perp}\subset\Bbbk G$, we can express 
$$e_{{\lambda}_{X^{\perp}}}=\sum_{\mu\in Y}e_{\mu}$$
in terms of the basis $\{e_{\mu}\}$ for $\Bbbk G$, indexed over some set $Y\subset\G$. Multiplying this expression by $h\in X^{\perp}$ we conclude that
$$\lambda(h)e_{{\lambda}_{X^{\perp}}}=\sum_{\mu\in Y}\mu(h)e_{\mu}.$$
This implies that $\mu_{X^{\perp}}=\lambda_{X^{\perp}}$, i.e. $\mu\in\lambda X$. Therefore the number of $\mu$'s occuring in the sum equals the dimension of $e_{{\lambda}_{X^{\perp}}}\Bbbk G$, which is $|G/X^{\perp}|=|X|$. It follows that the sum runs over all of $\lambda X$, as desired.\qed

%%%%%%%%%%%%%%%%%%%%%%%%%%%%%%%%%%%%%%%%%%%%%%%%%%%%%%%%%%%%%%%%%%%%%%%%%%%%%%%%%%%%%%%%%%%%%%%%%%%%%%%%%%%%%%%%%%%%%%%%%%%%%%%%%%%%%%%%%%%%%%%%%%%%%%%%
\section{Representations of Unlinked Liftings}\label{unlinked}
\subsection{Nilpotent unlinked liftings}\label{nilunliked}
\begin{Thm}\label{nilunlink}Let $H$ be an unlinked lifting and suppose $\epsilon_1=0=\epsilon_2$. Then
\begin{itemize}
\item[(a)] Every indecomposable projective module has the form $P_{\lambda}=H\e$ for some $\lambda\in\G$
\item[(b)] The central idempotents of $H$ are the $\e$ indexed by the cosets $\lambda X\in{\G}/X$.
\item[(c)] The Gabriel quiver for each block $H\ex$ has vertices corresponding to elements of $\lambda X$ and a pair of arrows
\begin{align*}\mu&\rightarrow\mu\chi_1\\
              \mu&\rightarrow\mu\chi_2\end{align*}
for every $\mu\in\lambda X$.
\item[(d)] $H$ is of wild representation type, unless $n_1=n_2=2$ in which case it is of tame representation type.
\end{itemize}
\end{Thm}
\pf Set $J=xH+yH$. By the hypothesis, $J$ is a nilpotent ideal, and so is the Jacobson radical of $H$. Therefore $H$ has $|G|$ indecomposable projective modules. Thus the projective modules $H\e,\;\lambda\in\G$ account for all of them. This proves (a).

We let $L_{\lambda}=H\e/J\e$ denote the corresponding simple $H$- module. It is easy to see that for all $r\in\Bbb N,\;\lambda\in\G$,
$$ J^r\e=\sum_{i+j\ge r}\k x^iy^j\e,$$
whence each $L_{\lambda}$ is one-dimensional, and $J\e/J^2\e$ is two-dimensional, being spanned by the images of $x$ and $y$. Note that the respective simple modules have weights $\lambda\chi_1$ and $\lambda\chi_2$. This says that $\text{Ext}_H^1(L_{\lambda},\, L_{\lambda\chi_i}),\,i=1,2$ are one - dimensional. As $\lambda\in\G$ is arbitrary, this accounts for the arrows labelled as in the statement. Thus the blocks are precisely indexed by the cosets as claimed. This completes the proof of (c).

Since $H$ is a basic algebra, the central idempotents are exactly sum of the primitive idempotents corresponding to the vertices in each connected component of the quiver. This yields (b).

$H$ is a skew group ring over the subalgebra generated by $x$ and $y$. One can easily see that we can replace $x$ by $a^{-1}x$ and assume that $x$ and $y$ commute. Thus $H$ is a skew group ring over a truncated polynomial ring, say $A$. Since $|G|$ is invertible in $\k$, the arguments  \cite[6.3]{ARS} show that $A$ and $H=AG$ have the same representation type. The assertion for the truncated polynomial rings follows from \cite{Dr}, whence (d).\qed 
              
%%%%%%%%%%%%%%%%%%%%%%%%%%%%%%%%%%%%%%%%%%%%%%%%%%%%%%%%%%%%%%%%%%%%%%%%%%%%%%%%%%%%%%%%%%%%%%%%%%%%%%%%%%%%%%%%%%%%%%%%%%%%%%%%%%%%%%%%%%%%%%%%%%%%%%%%
\subsection{Seminilpotent unlinked liftings}\label{semnilunlinked}
Next we describe indecomposable projective modules for unlinked data with exactly one $\epsilon_i=0$. Let $X_i$ denote the subgroup of $\G$ generated by $\chi_i$, and put $X=X_1X_2$ as before. We introduce another subgroup of $\G$ in this setting as follows. Let $N$ denote the subgroup $\langle a^{n_1}\rangle^{\perp}$ of $\G$. Note that $X\subset N$. This can be argued as follows. Since $\chi_1^{n_1}=\epsilon$ and $\chi_2(a)=\chi_1(b^{-1})$, we obtain $\chi_2(a^{n_1})=1$, which proves the inclusion.

We let $N/X_1$ denote the cosets modulo $X_1$ represented by elements of $N$ as usual, and write $(\G\setminus N)/X_1$ for the complementary set of cosets. Note that the disjoint sets $N/X_1$ and $(\G\setminus N)/X_1$ are stable under the action by multiplication by elements of $X$.
\begin{Thm}\label{seminilpunlinked} Suppose $H$ is an unlinked seminilpotent lifting with $\epsilon_1=1$ and $\epsilon_2=0$. Then
\begin{itemize}
\item[(a)] Let $\lambda\in N$. Then $He_{\lambda}$ is the projective cover of a one- dimensional $H$- module.
\item[(b)] Let $\mu\in\G\setminus N$. Then $He_{\mu}$ is the projective cover of a simple $n_1$- dimensional $H$- module.
\item[(c)] The isotypic component of $He_{\mu}$ is given by equivalence modulo $X_1$, consisting of the projectives $He_{\lambda},\lambda\in\mu X_1$.
\item[(d)] The $e_{\lambda},\lambda\in\G$ are a complete set of primitive orthogonal idempotents in $H$
\item[(e)] Let $\lambda\in N$. Then $H\ex$ is a block. The Gabriel quiver of this block has vertices corresponding to elements of $\lambda X$ and a pair of arrows
\begin{align*}\mu&\to\mu\chi_1\\
\mu&\to\mu\chi_2\end{align*}
for every $\mu\in\lambda X$ (with doubled arrows if $\chi_1=\chi_2$). The block $H\ex$ is of wild representation type, unless $n_1=n_2=2$, in which case it is of tame representation type.
\item[(f)] Let $\mu\in\G\setminus N$. Then $He_{\mu X}$ is a block. The Gabriel quiver of this block is cyclic and has vertices corresponding to cosets in $\mu X/X_1$ and arrows
$$\mu X_1\to\chi_2\mu X_1$$
corresponding to multiplication by $\chi_2$. The block $He_{\mu X}$ is a Nakayama algebra.
\item[(g)] The $\ex$ indexed by the cosets $\lambda X\in\G/X$ are a complete set of block idempotents of $H$.
\end{itemize}
\end{Thm}
\pf As noted above the subgroup $N$ contains $X$, so that the equivalence classes $\lambda X,\lambda\in N$ form a partition of $N$; taking $\lambda\in\G\setminus N$ instead, we similarly see that the classes $\text{mod}\,X$ form a partition of $\G\setminus N$. By Lemma \ref{centralid} the idempotents $\ex$ are all central. The two partitions combine to form a partition of $\G$ and we let 
\begin{align*}H_N&=\sum_{\lambda\in N}He_{\lambda}\\ H_{\G\setminus N}&=\sum_{\lambda\in\G\setminus N}He_{\lambda}\end{align*}
denote the complementary two-sided summands of $H$. Further let $J=xH_N+yH_N+yH_{\G\setminus N}$, which is a two-sided ideal of $H$ since $H_N$ and $H_{\G\setminus N}$ can be written as the sum of class subalgebras, each of the form $H\ex$.

First suppose that $\lambda\in N$. Then $x^{n_1}e_{\lambda}=0$ as $\lambda(a^{n_1})=1$. We also have $y^{n_2}=0$, hence the image $J$ in $H_N$ is nilpotent. Thus $J$ lies in the Jacobson radical of $H_N$. Next consider the projective modules $H\e$. It is plain that $J\e=\oplus_{i+j>0}\Bbbk x^iy^j\e$, so $H\e/J\e$ is one-dimensional, spanned by a vector of weight $\lambda$. This finishes the proof of (a).

Secondly, suppose that $\mu\in\G/N$ and consider the factor $He_{\mu}/Je_{\mu}$. Since $y^{n_2}=0$, it is clear that the image of $J$ in $H_{\G\setminus N}$ is nilpotent. In this case  $x^{n_1}e_{\mu}$ is a nonzero scalar and 
$$He_{\mu}/Je_{\mu}=\oplus_{0\le i\le n_1-1}\Bbbk x^ie_{\mu}+Je_{\mu}$$
Observe that the action of $x$ cyclically permutes the basis $\{x^ie_{\mu}+Je_{\mu}\}$, the vector $x^ie_{\mu}$ has weight $\mu\chi^i$ for all $i=0,\ldots, n_1-1$ and that these weights are pairwise distinct. It now follows easily that $He_{\mu}/Je_{\mu}$ is a simple $H$- module. For future reference let us denote this simple module by $L_{\mu}$. This finishes the proof of (b).

We show next that the isotypic component of $L_{\mu}$ consists of the simple modules $\{L_{\mu\chi}|\chi\in X_1\}$. Since the multiplicity of $L_{\mu}$ equals $\dim L_{\mu}=n_1$, we need to find $n_1$ modules isomorphic to $L_{\mu}$. If $\mu=\mu'\chi_1^j$, then the $H$- module map specified by $e_{\mu}+Je_{\mu}\mapsto x^je_{\mu'}+Je_{\mu'}$ gives an isomorphism $L_{\mu}\to L_{\mu'}$. This accounts for the isotypic component consisting of the $n_1$ mutually isomorphic projective (or simple) modules, which demonstrates (c).

To prove (d), first observe that the $e_{\lambda},\lambda\in\G$ are a complete set of orthogonal idempotents in $\Bbbk G$ by construction (\S \ref{abgroups}). The proofs of parts (a) and (b) show the projective modules of $H\e$ are all simple modulo the nilpotent ideal $J$. Thus $J$ is the radical of $H$, so the $\e$ are all primitive as well.

The proof of (e) is similar to the nilpotent case addressed in the previous theorem, and will be omitted.

We prove (f). Observe that
$$Je_{\mu}/J^2e_{\mu}=\oplus_{0\le i<n_1}\Bbbk x^iye_{\mu}+J^2e_{\mu}.$$
This factor is spanned by vectors $x^iye_{\mu}+J^2e_{\mu}$ with weights  $\mu\chi_2\chi^i_1$, which are permuted by the action of $x$, and annihilated by $y$. Thus it is apparent that $Je_{\mu}/J^2e_{\mu}$ is isomorphic to $L_{\mu\chi_2}$. This shows that the Gabriel quiver is as asserted in (f) and also that the block containing $L_{\mu}$ consists of the simples indexed by $\mu X_1X_2=\mu X$. Similar observations show that $J^ie_{\mu}/J^{i+1}e_{\mu}$ is isomorphic to $L_{\mu\chi_2^i}$  for $i=0,\ldots,n_2-1$. Therefore $He_{\mu}$ is uniserial for all $\mu\in\G/N$ and the proof of (f) is complete.

Finally, summing over the idempotents occuring in each block in the cases (e) and (f) (including multiplicities) produces every block idempotent as assserted in (g). This completes the proof of the theorem.\qed

%%%%%%%%%%%%%%%%%%%%%%%%%%%%%%%%%%%%%%%%%%%%%%%%%%%%%%%%%%%%%%%%%%%%%%%%%%%%%%%%%%%%%%%%%%%%%%%%%%%%%%%%%%%%%%%%%%%%%%%%%%%%%%%%%%%%%%%%%%%%%%%%%%%%%
\subsection{Unipotent unlinked liftings}

We now pass to the remaining unipotent case with $\epsilon_i=1$ for all $i$. In this case we require $\lambda(a^{n_1})\not=1\not=\lambda(b^{n_2})$. It is easy to see that such $\lambda$ exist. Recall $m=|X|$ and note that $m$ divides $n_1n_2$. Let $r={n_1n_2}/m$.
\begin{Thm}\label{unipotunlink} Let $H$ be an unlinked lifting and suppose $\epsilon_1=1=\epsilon_2$ with $\lambda(a^{n_1})\not=1\not=\lambda(b^{n_2})$. Then $H\ex$ is a semisimple algebra isomorphic to the direct sum of $r$ copies of $M_m(\k)$.
\end{Thm}
\pf Replace $y$ by $yb^{-1}$ so by abuse of notation we have $xy=yx$. We also abuse notation and write the images of generators in $H\ex$ as the generators $x,y$ themselves.

Since $x^{n_1}$ and $y^{n_2}$ are nonzero scalars, we can rescale and assume $x^{n_1}=y^{n_2}=1$. Let $U$ denote the subgroup of $H\ex$ generated by $x$ and $y$. For each $\alpha\in\widehat{U}$, let $f_{\alpha}$ be the primitive idempotent of $\k U$ associated to $\alpha$. By Proposition \ref{basis} $H\ex$ has basis
$$\{e_{\lambda\chi}x^{j_1}y^{j_2}|0\le j_i\le n_i,\chi\in X\},$$
so it follows that $H\e f_{\alpha}$ has basis
$$\{e_{\lambda\chi}f_{\alpha}|\chi\in X\}$$
We claim that $H\ex f_{\alpha}$ is a minimal left ideal. Let $\mu\in\lambda X$, and observe that the basis element $e_{\mu}f_{\alpha}$ has weight $\mu$. Furthermore,
\begin{align*}
x^{j_1}y^{j_2}e_{\mu}f_{\alpha}&=e_{\mu\chi}x^{j_1}y^{j_2}f_{\alpha}\\
                                &=\alpha(x^{j_1}y^{j_2})e_{\mu\chi}f_{\alpha}.
\end{align*}                                  
with $\chi=\chi_1^{j_1}\chi_2^{j_2}$. Therefore the given basis is permuted transitively by left multiplication by $U$ and the basis elements have distinct $G$-weights. These facts imply that each $H\ex f_{\alpha}$ is a minimal left ideal.

The decomposition 
$$H\ex=\bigoplus_{\alpha\in U}H\ex f_{\alpha}$$
demonstrates the semisimplicity of $H\ex$. Since $\dim H\ex=n_1n_2m$ and $\dim H\ex f_{\alpha}=m$, the proof is complete.\qed

We can describe isotypic components for the unlinked unipotent lifting next. Let $U$ be the subgroup generated by $x,y$ as in the previous result, and let $W$ be the kernel of the map
$$\pi:U\to X$$
given by $x\mapsto \chi_1,\,y\mapsto \chi_2$. For $\alpha\in\widehat U$, let $\alpha_W\in\widehat W$ denote the restriction of $\alpha$ to $W$.
\begin{Thm}\label{unipotent}Let $\alpha,\beta\in\widehat W$. Then $H\ex f_{\alpha}\cong H\ex f_{\beta}$ if and only if $\alpha_W=\beta_W$.
\end{Thm}
\pf To prove this result it suffices to show that $f_{\beta}H\ex f_{\alpha}=0$ if and only if $\alpha_W\not=\beta_W$. To this end we calculate
\begin{align}
f_{\beta}e_{\mu}f_{\alpha}&=\frac{1}{n_1n_2}\sum_{z\in U}\beta^{-1}(z)ze_{\mu}f_{\alpha}\nonumber\\
                          &=\frac{1}{n_1n_2}\sum_{z\in U}(\beta^{-1}\alpha)(z)e_{\mu\pi(z)}f_{\alpha}\label{intertwine}
                          \end{align}
for all $\mu\in\lambda X$, where the second expression is obtained by noting that $ze_{\mu}f_{\alpha}=e_{\mu\pi(z)}zf_{\alpha}=\alpha(z)e_{\mu\pi(z)}f_{\alpha}$.

Put $\omega=\beta^{-1}\alpha$ and assume that $\alpha_W\not= \epsilon$. Let $T$ denote a transversal for $W$ in $U$. Now we can write (\ref{intertwine}) in the form 
$$\frac{1}{n_1n_2}\sum_{t\in T}\omega(t)(\sum_{w\in W}\omega(w))e_{\mu t}f_{\alpha}$$
which is zero by \eqref{ort}.

Moving to the opposite inclusion, suppose that $\omega_W=\epsilon$. Now (\ref{intertwine}) can be written 
\begin{equation}\label{i2}\frac{W}{n_1n_2}\sum_{t\in T}\omega(t)e_{\mu\pi(t)}f_{\alpha}=\frac{|W|}{n_1n_2|G|}\sum_{g\in G}\mu(g^{-1})[\sum_{t\in T}\omega(t)\pi(t)(g^{-1})]gf_{\alpha}\end{equation}
using the definition of the $\e$. By duality for finite abelian groups $G/{X^{\perp}}\simeq\widehat X$ under the evaluation map $g\mapsto \text{ev}_g:\chi\mapsto \chi(g),\,\chi\in X$. Since $U/W\simeq X$ there exists $g_0\in G$ unique modulo $X^{\perp}$ such that $\omega(t)=\pi(t)(g_0)$ for all $t\in T$. Therefore the inner sum in \eqref{i2} becomes
\begin{equation}\label{i3}\sum_{t\in T}\pi(t)(g_0g^{-1})\end {equation}
As $t$ runs over $T$, $\pi(t)$ runs over $X$. Thus the sum in \eqref{i3} is in fact
$$\sum_{\phi\in X}\phi(g_0g^{-1})=\begin{cases}0, &g_0g^{-1}\notin X^{\perp}\\|X|,&g_0g^{-1}\in X^{\perp}\end{cases}$$
 by using \eqref{ort} and duality $G/{X^{\perp}}\simeq\widehat X$ mentioned above. Thus we can harmlessly assume that $g$ runs over $g_0X^{\perp}$ and then \eqref{i2} reduces to 
$$\frac{1}{|G|}\sum_{h\in X^{\perp}}\mu^{-1}(g_0h)g_0hf_{\alpha}=\frac{|X^{\perp}|}{|G|}g_0\mu^{-1}(g_0)\ex f_{\alpha},$$
the second equality by Lemma \ref{idempotents}. This expression is clearly nonzero, so the proof is complete.\qed 
%%%%%%%%%%%%%%%%%%%%%%%%%%%%%%%%%%%%%%%%%%%%%%%%%%%%%%%%%%%%%%%%%%%%%%%%%%%%%%%%%%%%%%%%%%%%%%%%%%%%%%%%%%%%%%%%%%%%%%%%%%%%%%%%%%%%%%%%%%%%%%%%%%%%%%
%%%%%%%%%%%%%%%%%%%%%%%%%%%%%%%%%%%%%%%%%%%%%%%%%%%%%%%%%%%%%%%%%%%%%%%%%%%%%%%%%%%%%%%%%%%%%%%%%%%%%%%%%%%%%%%%%%%%%%%%%%%%%%%%%%%%%%%%%%%%%%%%%%%%%%

\section{Representations of Linked Liftings}\label{linked}
\subsection{Linked liftings}\label{LL}
Let $H$ be  a linked lifting with datum $\D$. Setting $q=\chi_2(a), \chi=\chi_1=\chi_2^{-1},\,\text{and}\, n=|q|$ one can see readily that $\chi(a)=\chi(b)=q^{-1}$ and $n_1=n_2$. Modifying the defining relations in \ref{Types} we arrive at a presentation of $H=H(\D)$ by $G$ together with the relations
\begin{align}x^n&=\epsilon_1(a^{n}-1)\label{xpower}\\
              y^n&=\epsilon_2(b^{n}-1)\label{ypower}\\
              gx&=\chi(g)xg\label{char1}\\
              gy&=\chi^{-1}(g)yg\label{char2}\\
              xy -qyx&=\gamma(ab-1)\label{qcom}\end{align}
with the coalgebra structure as given in \ref{Types}.              

A very important feature of linked liftings is the presence of a special central element analogous to the Casimir element of the classical Lie theory.
\begin{Lem}\label{cas}\begin{align} C&=(a^{-1}x)y-\frac{\gamma}{q-1}(a^{-1}+qb)\label{cas1}\\
                                     &=y(a^{-1}x)-\frac{\gamma}{q-1}(b+qa^{-1})\label{cas2}.\end{align}
                                     is a central element of $H$

\end{Lem}
\pf The second formula follows from the first by \eqref{qcom}. The first claim of the lemma is a straightforward verification. We sketch details. First, by \eqref{char1} and \eqref{char2} $xy$ commutes with every g$\in G$. Further,
\begin{align*}a^{-1}x(yx)&=xqa^{-1}(yx)\\&=x(qa^{-1}(q^{-1}xy-q^{-1}\gamma(ab-1))\\&=x(a^{-1}xy-\gamma(b-a^{-1}).\end{align*}
by \eqref{char1} and \eqref{qcom}, and
$$(a^{-1}+qb)x=x(qa^{-1}+b).$$
by \eqref{char1}. It follows that
\begin{align*}Cx&=x(a^{-1}xy-\gamma b(1+\frac{1}{q-1})-\gamma a^{-1}(\frac{q}{q-1}-1))\\&=xC.\end{align*}
The equality $yC=Cy$ is established similarly using the second expression for $C$.\qed

We proceed to a very useful property of linked liftings.
\begin{Prop}\label{symalg}Every linked lifting of a quantum plane is a symmetric algebra.
\end{Prop}
\pf Let $H$ be a linked lifting of the quantum plane. By a well-known result in \cite{OS} the assertion is equivalent to $H$ being unimodular with  inner square of the antipode.

We first show that $S^2$ is inner. It is easy to see that $S^2(x)=-a^{-1}xa$ and $S^2(y)=-b^{-1}yb$. Since $ab^{-1}\in\ker\chi,ab^{-1}$ is central in $H$. Therefore $S^2(y)=a^{-1}ya$ as well. It follows immediately that $S^2(h)=-a^{-1}ha$ for all $h\in H$.

It remains to show that $H$ is unimodular. To do this we provide an integral $I:=ex^{n-1}y^{n-1}$ where $e=|G|^{-1}\sum_{g\in G}g$.

Let $g\in G$. Since $\chi_1\chi_2=\epsilon$, it is immediate that $x^{n-1}y^{n-1}$ commutes with $g$. Hence 
$$gI=Ig=\epsilon(g)I.$$
We show that $xI=Iy=0$ next. This is clear if $x^n=0=y^n$, so assume, say, $x^n=a^n-1$ with $\chi^n=\epsilon$. Then $\chi(a^n)=1$ and hence $a^n$ is central in $H$. Therefore
\begin{align*}xI&=x(x^{n-1}y^{n-1}e)\\&=x^ny^{n-1}e\\&=(a^n-1)y^{n-1}e\\&=0.\end{align*}
Similarly $Iy=0$.

It remains to show that $Ix=0=yI$. We shall prove $yI=0$. 

Note that as $x^{n-1}y^{n-1}$ has weight $\epsilon$, $[y,x^{n-1}y^{n-1}]_c=yx^{n-1}y^{n-1}-x^{n-1}y^n$. Therefore
$$yI=x^{n-1}y^ne+[y,x^{n-1}y^{n-1}]_ce.$$
Since $y^n=\epsilon_2(b^n-1)$ we have $y^ne=0$. Furthermore, using Proposition \ref{productrule} we have for the second term
\begin{align*}[y,x^{n-1}y^{n-1}]_c&=[y,x^{n-1}]_cy^{n-1}+\chi^{n-1}(b)x^{n-1}[y,y^{n-1}]_c\\
&=\gamma q^{-(n-1)}(n-1)_qx^{n-2}(q^{-n+2}ab-1)y^{n-1}\\&\quad + \chi^{n-1}(b)x^{n-1}(1-q^{-(n-1)})y^n\end{align*}
Since $aby^{n-1}=q^{2(n-1)}y^{n-1}ab$ and $y^ne=0$ we conclude that for some $h\in H$
$$yI=h(ab-1)e$$
which is zero, completing the proof of the Theorem.\qed

We proceed to a description of the class subalgebras of $H$. Let $X$ be the subgroup of $G$ generated by $\chi$ and put $N=|X|$. Since $q=\chi^{-1}(a)$, $n|N$. We set $m=N/n$.
\begin{Thm}\label{linkclsalg}Let $H$ be  a linked lifting and $\lambda\in\G$.
\begin{itemize}
\item[(a)] The algebra $H\ex$ is generated by elements $E,F,K$ subject to the relations
\begin{align}E^n&=\pm\epsilon_1\lambda(a^{-n}-1)\label{Epower}\\
             F^n&=\epsilon_2\lambda(b^n-1)\label{Fpower}\\
             K^N&=\kappa\\
             EF&-FE=\eta(K^{-m}-K^m)\label{comrel}\\
             KE&=\theta EK\label{KErel}\\
             KF&=\theta^{-1}FK\label{KFrel}\end{align}
 where $\theta$ is a primitive $N^{\text th}$ root of one, $\kappa,\eta$ are nonzero scalars and $\epsilon_i=0, 1$.
 \item[(b)] The algebra $H\ex$ is of dimension $Nn^2$.
 \end{itemize}\end{Thm}
 \pf We begin by modifying the defining relations of $H$. Let $x'=a^{-1}x$. Then \eqref{qcom} becomes
 \begin{align}x'y-yx'=\gamma(b-a^{-1})\label{commutator}.\end{align}
 A simple calculation gives $(ax')^n=q^{\binom n2}a^nx'^n=\pm a^nx'^n$. By rescaling $x'$ if necessary we may assume that \eqref{xpower} holds for $x'$ in place of $x$.
 
Let $E=\ex x'$ and $F=\ex y$. Multiplying \eqref{xpower} and \eqref{ypower} by $\ex$ and noting that $a^n,b^n\in X^{\perp}$, we obtain \eqref{Epower} and \eqref{Fpower}. Further, note that $\chi(G)$ is finite, hence cyclic subgroup of $\Bbbk^{\bullet}$ and that $q\in\chi(G)$ has order $n$. So we can choose $\theta\in\Bbbk$ to be a generator of $\chi(G)$ and $g\in G$ such that 
\begin{align*}\chi(g)&=\theta\\\theta^m&=q.\end{align*}
Since $\chi(G)=G/{X^{\perp}},\theta$ has order $N$. Also $a^{-1}=h_1g^m$ and $b=h_2g^{-m}$ for some $h_i\in X^{\perp}$.

Let $\overline{g}=\ex g, \alpha=\lambda(h_1)$ and $\beta=\lambda(h_2)$. Multiplying \eqref{commutator} by $\ex$ and using Lemma \ref{idempotents} we obtain
$$EF-FE=\gamma(\beta\overline{g}^{-m}-\alpha\overline{g}^m).$$
Pick $\mu\in\Bbbk$ such that $\mu=\alpha/\beta$ and set $K=\mu\overline{g}$. One can check that
$$\beta\overline{g}^{-m}-\alpha\overline{g}^m=\sqrt{\alpha\beta}(K^{-m}-K^m).$$
Setting $\eta=\gamma\sqrt{\alpha\beta}$ results in \eqref{comrel}. This completes the proof of (a).

It is clear that any algebra so presented has dimension at most $Nn^2$. In view of the fact that $[G:X^{\perp}]=|\widehat{X}|=N$ and Proposition \ref{basis} we see that the dimension of $H\ex$ is precisely this upper bound.\qed
 
%%%%%%%%%%%%%%%%%%%%%%%%%%%%%%%%%%%%%%%%%%%%%%%%%%%%%%%%%%%%%%%%%%%%%%%%%%%%%%%%%%%%%%%%%%%%%%%%%%%%%%%%%%%%%%%%%%%%%%%%%%%%%%%%%%%%%%%%%%%%%%%%%%%%%%%%                
\subsection{Casimir element}\label{casimir}
Write $C_{\lambda}$ for the image $C\ex$ of $C$ in $H\ex$. In this subsection we determine the minimal polynomial for each $C_{\lambda}$. First, we set up some notation. Explicitly we have
\begin{align}C_{\lambda}&=EF-\frac{\eta}{q-1}(K^m+qK^{-m})\label{cas1'}\\
                        &=FE-\frac{\eta}{q-1}(K^{-m}+qK^m)\label{cas2'}.\end{align}
Let ${\eta}'=\frac{\eta}{q-1}$ and $D=K^{-m}+qK^m$. Denote by $A_{\lambda}$ the subalgebra of $H\ex$ generated by $K$ and define the $\k$- automorphism of $A_{\lambda}$ by the rule $\tau(K^i)={\theta}^i K^i$. Notice that the action of $\tau$ extends to an automorphism of $A_{\lambda}[t]$ with fixed ring $\k [t]$ and $aE=E\tau(a)$ for all  $a\in A_{\lambda}$. Further, set
\begin{equation}\label{roots}R_{{\kappa},N}:=\{\alpha\in\k|\alpha^N=\kappa\}\;\text{and}\;R_N=R_{1,N}.\end{equation}
For every $\rho\in R_{\kappa,N}$, let $\phi_{\rho}:A_{\lambda}\to\k$ denote the evaluation homomorphism given by sending $K$ to $\rho$.               
\begin{Lem}\label{minpol}For every $\rho\in R$
$$f_{\lambda}(t)=\prod_{i=0}^{n-1}(t+\eta'\phi_{\theta^i\rho}(D))-E^nF^n$$
is the minimal polynomial of $C_{\lambda}$ over $\k$.
\end{Lem}
\pf We first note that the minimal polynomial of $C_{\lambda}$ is of degree $n$. For otherwise, suppose $g(C_{\lambda})=0$ for some monic $g(t)\in\k[t]$ with $\deg g(t)=m<n$. By \eqref{comrel} and Lemma \ref{commutators}
$$(EF)^m=E^mF^m+\sum_{1\le i<m}a_iE^iF^i$$
for some $a_i\in A_{\lambda}$. Hence $C_{\lambda}^m=E^mF^m+\sum_{i<m}b_iE^iF^i$ for some $b_i\in A_{\lambda}$ and we deduce that
$$E^mF^m+\sum_{i<m}c_iE^iF^i=0$$
for some $c_i\in A_{\lambda}$. But this contradicts the fact that $\{E^iF^jK^k|0\le i,j,k<n\}$ is a basis for $H\ex$.

We show next that $f_{\lambda}(C_{\lambda})=0$. Notice that from the definition of $\tau$ we have $aE^i=E^i\tau^i(a)$ for all $a\in A_{\lambda}$. Also, by the definition of $C_{\lambda}$, $FE=C_{\lambda}+\eta'D$. More generally for all $i>0$,
\begin{align*}F^iE^i&=F^{i-1}(FE)E^{i-1}\\
                    &=F^{i-1}E^{i-1}(C_{\lambda}+\eta'\tau^{i-1}(D)).\end{align*}
Iterating we see that
$$F^nE^n=(C_{\lambda}+\eta'D)(C_{\lambda}+\eta'\tau(D))\cdots(C_{\lambda}+\eta'\tau^{n-1}(D)).$$
This says precisely that $C_{\lambda}$ satisfies $f_{\lambda}(t)$.

It remains to show that $f_{\lambda}(t)\in\k[t]$. Write $f_{\lambda}(t)=\sum_{j=0}^ne_jt^{n-j}-E^nF^n$ where $e_j$ is the $j$th elementary symmetric function of $\eta'\tau^i(D)$. As $\tau$ has finite order, the $\eta'\tau^i(D)$ are permuted by $\tau$, hence the $e_j$ are fixed by it, and so $e_j\in\k$.

On the other hand a direct inspection using $\theta^m=q$ gives $\phi_{\rho}(\tau^i(D))=\phi_{\theta^i\rho}(D)$. Since $f_{\lambda}(t)\in\k[t]$ and $E^nF^n\in\k$ we compute
\begin{align*}f_{\lambda}=\phi_{\rho}(f_{\lambda}(t))&=\prod_{i=0}^{n-1}(t+\eta'\phi_{\rho}(\tau^i(D))-E^nF^n\\
                                                     &=\prod_{i=0}^{n-1}(t+\eta'\phi_{\theta^i\rho}(D))-E^nF^n\end{align*}.\qed
                    
%%%%%%%%%%%%%%%%%%%%%%%%%%%%%%%%%%%%%%%%%%%%%%%%%%%%%%%%%%%%%%%%%%%%%%%%%%%%%%%%%%%%%%%%%%%%%%%%%%%%%%%%%%%%%%%%%%%%%%%%%%%%%%%%%%%%%%%%%%%%%%%%%%%%%%%%
\subsection{The unipotent case}\label{linkpotent} 
\begin{Thm}\label{potentclsalg} Suppose $E^nF^n\neq 0$. Then $H\ex$ is isomorphic to the algebra of $n\times n$ matrices over $\k[C_{\lambda}]$.
\end{Thm}
\pf Since $E^n$ and $F^n$ are both nonzero by hypothesis, we have $\epsilon_1=\epsilon_2=1$ in the datum for $H$. It follows that $\chi^n=\epsilon$ and hence $N=n$. Let $B^-=\k [K,F]$ denote the algebra generated by $K$ and $F$. It is easy to see that $B^-$ is isomorphic to the skew group algebra where the group generated by $K$ acts faithfully on $\k[F]$ by multiplication by powers of $q^{-1}$. By a standard argument e.g. \cite{Rad}, we obtain $B^-\cong M_n(\k)$.

Since $F$ is invertible, we have $E=F^{-1}C_{\lambda}+\eta'F^{-1}D$. Therefore $H\ex$ is generated by $K, F,C_{\lambda}$. Since $\dim H\ex=n^3$, we see that $$\{K^iF^jC_{\lambda}^k|0\le i,j,k<n\}$$ is a basis for $H\ex$. It follows that $H\ex\cong B^-\otimes \k[C_{\lambda}]\cong M_n(\k[ C_{\lambda}])$.\qed
\begin{Cor} Suppose $E^nF^n\neq 0$. Then every simple $H\ex$-module is $n$-dimensional.\qed
\end{Cor}
%%%%%%%%%%%%%%%%%%%%%%%%%%%%%%%%%%%%%%%%%%%%%%%%%%%%%%%%%%%%%%%%%%%%%%%%%%%%%%%%%%%%%%%%%%%%%%%%%%%%%%%%%%%%%%%%%%%%%%%%%%%%%%%%%%%%%%%%%%%%%%%%%%%%%%%

\subsection{The semi-nilpotent case}\label{linkseminilp}
We now consider the case $F^n\neq 0$ and $E^n=0$ with varying assumptions on $\kappa$. As in the previous theorem we have $n=N$. In the following theorem we put 
$$V=\k[v]/{(v^2)}$$
for the truncated polynomial algebra. For an algebra $A$ we let $A^j$ denote the direct sum of $j$ copies of the $A$.
\begin{Thm}\label{lsn} Suppose $F^n\neq 0$ and $E^n=0$.
\begin{itemize}
\item[(a)] If $\kappa^2\neq 1$, then $H\ex$ is isomorphic to $M_n(\k)^n$.

Assume now that $\kappa^2=1$
\item[(b)] If $n$ is odd , then $H\ex$ is isomorphic to
$$M_n(\k)\oplus M_n(V)^{\frac{n-1}{2}}.$$

Assume further that $n$ is even
\item[(c)] If $\kappa=1$, then $H\ex$ is isomorphic to 
$$M_n(V)^{\frac{n}{2}}.$$
\item[(d)] If $\kappa=-1$, then $H\ex$ is isomorphic to 
$$M_n(\k)^2\oplus M_n(V)^{\frac{n-2}{2}}.$$
\end{itemize}
\end{Thm}
\pf The argument used in the preceding Theorem shows that $H\ex\cong M_n(\k[C_{\lambda}])$. Let $\rho\in R_{\kappa,n}$. Then $R_{\kappa,n}=\{q^i\rho|i=0,1\ldots,n-1\}$. As $E^n=0$, we obtain the factorization
$$f_{\lambda}=\prod_{\rho\in R_{\kappa,n}}(t+\eta'\phi_{\rho}(D))=\prod_{\rho\in R_{\kappa,n}}(t+\eta'D(\rho)),$$
using the fact that $\theta=q$. So the asserted results in (a)- (d) are a question of determining the multiplicity of the roots $-\eta'(D(\rho))$. A simple computation gives that for $\rho_1,\rho_2\in R_{\kappa,n}$,
\begin{equation}\label{multiplicity} D(\rho_1)=D(\rho_2)\Leftrightarrow \rho_1\rho_2=q^{-1}.\end{equation}
In case $\kappa^2\neq 1$ this condition does not hold for any two roots, whence we get (a).

Suppose $\kappa^2=1$. Then \eqref{multiplicity} can be rephrased by stating that the root $-\eta'D(\rho)$ occurs with multiplicity one if and only if $\rho^2=q^{-1}$. Otherwise the root $-\eta'D(\rho)=-\eta'D(q^{-1}\rho^{-1})$ occurs twice. It follows that the primary decomposition of $\k [C_{\lambda}]$ is
\begin{equation}\label{primary}\k[C_{\lambda}]=\k^s\oplus(\oplus\k[t]/{(t+\eta'D(\rho))^2}\end{equation}
where $s$ is the number of simple roots. 

It remains to determine $s$. We note that for $\rho,\zeta\in R_{\kappa,n}$ $\zeta\rho^{-1}\in R_n$ and $\rho^2=\zeta^2$ if and only if $(\zeta\rho^{-1})^2=1$. The mapping $\sigma:R_n\to R_n,\;q^i\mapsto q^{2i}$ is an isomorphism if $n$ is odd, while $\text{im}\,\sigma=R_n^2$ for even $n$. If $\kappa=-1$, let $\pi$ be a primitive $2n^{\text{th}}$ root of $1$ such that $\pi^2=q$. Then the preceding paragraph shows that $s=1$ if n is odd, $s=0$ if $n$ is even and $\kappa=1$, and $s=2$ if $n$ is even and $\kappa=-1$.\qed
%%%%%%%%%%%%%%%%%%%%%%%%%%%%%%%%%%%%%%%%%%%%%%%%%%%%%%%%%%%%%%%%%%%%%%%%%%%%%%%%%%%%%%%%%%%%%%%%%%%%%%%%%%%%%%%%%%%%%%%%%%%%%%%%%%%%%%%%%%%%%%%%%%%%%%%%

\subsection{Standard cyclic modules}\label{babyverma}
We continue with the linked case in summands $H\ex$ with $E^n=0$. We explicitly construct some induced modules which turn out to be simple for the case where $F^n\neq 0$. In the next section we will assume $E^n=0=F^n$ and there these modules will play a crucial role. Let $B^+$ denote the subalgebra of $H\ex$ generated by $E,K$. For every $\rho\in R_{\kappa,N}$ we make $\k$ a $B^+$- module denoted by $\k_{\rho}$ by setting $K.1_{\rho}=\rho$ and $E.1_{\rho}=0$ where $1_{\rho}$ is identified with $1\in\k$. We define the $H\ex$- module $Z(\rho)$ by
$$Z(\rho)=H\ex\otimes_{B^+}\k_{\rho}$$
Since $H\ex$ is a free $B^+$- module with basis $\{F^i|0\le i\le n-1\}$, we have a $\k$- basis for $Z(\rho)$
$$\{F^i\otimes 1_{\rho}\}$$
and we shall write $w_i=F^i\otimes 1_{\rho}$. This is the {\em standard basis} for $Z(\rho)$. Using a conventional argument, it follows that each simple module is a homomorphic image of some $Z(\rho)$. When $F$ is invertible and the $F^i\otimes 1_{\rho}$ have distinct $K$- weights, it is immediate that $Z(\rho)$ is simple.

\begin{Prop}\label{simpleinduced}In case $E^n= 0,\;\{Z(\rho)|\rho\in R\}$ is a set of representatives of simple $H\ex$- modules. In addition, $Z(\rho)\cong Z(\zeta)$ if and only if $C_{\lambda}.1_{\rho}=C_{\lambda}.1_{\zeta}$.
\end{Prop}
\pf  For every $\rho\in R_{\kappa,n}$ we put $C_{\rho}$ for the summand of \ref{primary} corresponding to the root $-\eta'D(\rho)$ and let $\epsilon_{\rho}$ be the primitive idempotent generating $C_{\rho}$. Clearly $\epsilon_{\rho}I_n$ is the unity of $M_n(C_{\rho})$. Since $M_n(C_{\rho})$ has a unique isomorphism class of simple modules, $Z(\rho)$ is in the class if and only if the restriction $\epsilon_{\rho}|_{Z(\rho)}$ is the identity. Since $\epsilon_{\rho}$ is a polynomial in $C_{\lambda}$ taking the value $1$ at $-\eta'D(\rho)$ and zero at any other root of $f_{\lambda}(t)$, and as $C_{\lambda}\cdot 1_{\rho}=-\eta'D(\rho)$, the result follows.\qed 
 
If we assume that $F^n=0$ instead we can construct a $H\ex$- module in the obvious manner, inducing from the subalgebra $B^-$
and obtain $H\ex$ - modules which we denote by $Z'(\rho)$. These modules are simple if $E^n\neq 0$. We have a $\k$- basis for $Z'(\rho)$
$$\{E^i\otimes 1_{\rho}|0\le i\le n-1\}$$
and we shall write $v_i=E^i\otimes 1_{\rho}$. This is called the {\em standard basis} for $Z'(\rho)$.
%%%%%%%%%%%%%%%%%%%%%%%%%%%%%%%%%%%%%%%%%%%%%%%%%%%%%%%%%%%%%%%%%%%%%%%%%%%%%%%%%%%%%%%%%%%%%%%%%%%%%%%%%%%%%%%%%%%%%%%%%%%%%%%%%%%%%%%%%%%%%%%%%%%%%%%%

\subsection{The nilpotent case}\label{nilpotent}
From now on we assume $E^n=0=F^n$. We use the notation of the previous subsection, including the modules $Z'(\rho)$ and $Z(\rho)$ and their standard bases. Recall that $N=nm$ and the definition of $R_{\kappa,N}$. We will refer to elements of $R_{\kappa,N}$ as roots. 

\begin{Def} Let $p:R_{\kappa,N}\to R_{\kappa^2,n}$ be defined by $p(\rho)=\rho^{2m}$. We call $\rho\in R_{\kappa,N}$ {\em exceptional} if $\rho^{2m}=q^{n-1}$. For every $\rho\in R_{\kappa,N}$ we define integers $e(\rho)$ and $e'(\rho)$ by the rule
$$p(\rho)=q^{e(\rho)}=q^{-e'(\rho)}.$$
\end{Def}
\begin{Def} Let $M$ be an $H\ex$- module. We say that $0\neq v\in M$ has {\em weight} $\pi\in R$ if $K.v=\pi v$. A weight vector $v$ is said to be 
$E$-{\em trivial} if $E.v=0$.
\end{Def}
We preface a description of the induced modules with a commutation formula similar to one for quantized universal enveloping algebra for $\mathfrak{sl}_2$ \cite{KA}. The proof is entirely similar to that of Proposition \ref{productrule} (b) and will be omitted.
\begin{Lem}\label{commutators}For every $s\ge 1$
\begin{itemize}
\item[(1)] $[E,F^s]=\eta (s)_qF^{s-1}(K^{-m}-q^{-(s-1)}K^m)$\label{com1}.
\item[(2)] $[F,E^s]=\eta (s)_qE^{s-1}(K^m-q^{-(s-1)}K^{-m})$\label{com2}.
\end{itemize}
\end{Lem}
%%%%%%%%%%%%%%%%%%%%%%%%%%%%%%%%%%%%%%%%%%%%%%%%%%%%%%%%%%%%%%%%%%%%%%%%%%%%%%%%%%%%%%%%%%%%%%%%%%%%%%%%%%%%%%%%%%%%%%%%%%
\begin{comment}
\pf (1) The formula holds for $s=1$ by \eqref{comrel}. We induct on $s$ assuming the formula holds for $s$. We compute
\begin{align*}[E,F^{s+1}]&=F^s[E,F]+[E,F^s]F\\&=\eta F^s(K^{-m}-K^m + q^{1-s}(s)_q(q^{s-1}\theta^m K^{-m}-\theta^{-m}K^m)\\
\intertext {(using the induction hypothesis, \eqref{KFrel} and $\theta^m=q$)}
 &=\eta F^s(K^{-m}-K^m+(s)_qqK^{-m}-q^{-s}(s)_qK^m\\&=\eta q^{-s}(s+1)_q F^s(q^sK^{-m}-K^m)\end{align*}
the last equality by the identities $1+(s)_qq=(s+1)_q$ and $q^s+(s)_q=(s+1)_q$.

Part (2) is proven similarly using the relation \eqref{KErel}\qed
\end{comment}
%%%%%%%%%%%%%%%%%%%%%%%%%%%%%%%%%%%%%%%%%%%%%%%%%%%%%%%%%%%%%%%%%%%%%%%%%%%%%%%%%%%%%%%%%%%%%%%%%%%%%%%%%%%%%%%%%%%%%%%%%%%%%%%%%%%%%%%%%%%%

\begin{Prop}\label{reducible} \begin{itemize}\item [(a)] $Z'(\rho)$ is simple if and only if\newline $e'(\rho)= n-1$. 
If not, $Z'(\rho)$ has a unique nonzero proper submodule generated by $v_{e'+1}$.
\item[(b)] $Z(\rho)$ is simple if and only if $e(\rho)= n-1$. 
If not, $Z(\rho)$ has a unique nonzero proper submodule generated by $w_{e+1}$.\end{itemize}
\end{Prop}
\pf Recall the standard basis $v_i=E^i\otimes 1_{\rho},\,i=0,1,\ldots,n-1$. We show first that $Z'(\rho)$ has a nonzero proper submodule if and only if there is an $F$-trivial $v_s$ with $s>0$. 

Since $v_i$ has weight $\theta^i\rho$, the $v_i$ have distinct weights. Thus, a proper submodule $M$ of $Z'(\rho)$ is the span of the $v_i$ that lie in $M$. Let $s=\text{min}\{i|v_i\in M\}$. By Lemma \ref{commutators} (2) $F.v_s=0$, and $s>0$ because $M$ is proper. Conversely, if $F.v_s=0$ for some $s>0$, then $H\ex v_s$ is a nonzero proper submodule of $Z'(\rho)$. 

Using Lemma \ref{commutators}(2) we quickly deduce that (for $s>0$) $F.v_s=0$ if and only if $\rho^{2m}=q^{1-s}$. Setting $e'=s-1$, we arrive at the desired equivalence.

It follows from $\rho^{2m}=q^{1-s}$ that $H\ex v_s$ contains no $F$-trivial vectors other than $v_s$. Thus $H\ex v_s$ is the unique proper nonzero submodule of $Z'(\rho)$.

The proof of (b) is similar using Lemma \ref{commutators} (1).\qed

We dispose easily of the case where $\kappa^2\neq 1$.
\begin{Thm}\label{dumbcase}Assume that $\kappa^2\neq 1$. then $H\ex\cong M_n(\k)^N$.
\end{Thm}
\pf  Suppose that $Z(\rho)$ is not simple. Then the result above says that $\rho^{2m}=q^e$,which implies $\kappa^2=\rho^{2N}=\rho^{2mn}=q^{en}=1$. We conclude that every $Z(\rho)$ is simple. Furthermore, since each $Z(\rho)$ contains a unique $E$- trivial vector of weight $\rho$, we see that the $Z(\rho)$ are pairwise nonisomorphic. As $\dim H\ex=Nn^2$, the proof is complete.\qed
%%%%%%%%%%%%%%%%%%%%%%%%%%%%%%%%%%%%%%%%%%%%%%%%%%%%%%%%%%%%%%%%%%%%%%%%%%%%%%%%%%%%%%%%%%%%%%%%%%%%%%%%%%%%%%%%%%%%%%%%%%%%%%%%%%%%%%%%%%%%%%%%%%%%%%%
\subsection{ {\boldmath$\kappa^2=1$}}\label{mainnilcase}

Let $\text{rad}\, Z(\rho)$ (resp. $\text{rad}\, Z'(\rho)$) denote the proper (possibly zero) submodule of $Z(\rho)$ (resp. $Z'(\rho)$). Further let
\begin{align*}L(\rho)&=Z(\rho)/{\text{rad}\, Z(\rho)}.\\
              L'(\rho)&=Z'(\rho)/{\text{rad}\, Z'(\rho)}.\end{align*} 
\begin{Prop}\label{extensions}Assume $\kappa^2=1$.
\begin{itemize}
\item [(a)] The simple modules $L(\rho),\rho\in \R$ are a full set of representatives of simple $H\ex$- modules.
\item [(b)] The simple modules $L'(\rho),\rho\in \R$ are another set of representatives of simple $H\ex$- modules.
\item [(c)] If $Z(\rho)$ is not simple, $Z(\rho)$ is a nonsplit extension of $L(\theta^{-(e(\rho)+1)}\rho)$ by $L(\rho)$.
\item [(d)] If $Z'(\rho)$ is not simple, $Z'(\rho)$ is a nonsplit extension of $L(\theta^{n-1}\rho)$ by $L(\theta^{e'(\rho)}\rho)$.
\item [(e)] $\dim L(\rho)=e(\rho)+1$.
\item [(f)] $\dim L'(\rho)=e'(\rho)+1$.
\end{itemize}
\end{Prop}
\pf (a)-(b). Every simple $H\ex$- module contains an $E$- trivial vector $v$, say of weight $\rho\in\R$. The map given by $1_{\rho}\mapsto v$ induces a nontrivial homomorphism $Z(\rho)\to M$. Since $L(\rho)$ contains a unique line of $E$-trivial vectors of weight $\rho$ $,\;L(\rho)\ncong L(\rho')$ for $\rho\neq \rho'$. This proves (a). The proof of (b) is similar.

(c) and (e). By Proposition \ref{reducible} (b) $\text{rad}\,Z(\rho)$ is generated by an $E$- primitive $w_{e(\rho)+1}$, hence is the span of $\{w_i|e(\rho)+1\le i\le n-1\}$. This proves (e). Further, the weight of $w_{e(\rho)+1}$ is $\theta^{-(e(\rho)+1)}\rho$, so that $\text{rad}\,Z(\rho)\cong L(\theta^{-(e(\rho)+1)}\rho)$. This proves (c).

(d) and (f). By Proposition \ref{reducible} (a) $\text{rad}\,Z'(\rho)$ is the span of $\{v_i|e'(\rho)+1\le i\le n-1\}$. Now, up to a scalar multiple, the only $F$- trivial weight vector of $\text{rad}\,Z'(\rho)$ is $v_{n-1}$ and the only $F$-trivial weight vector in $L'(\rho)$ is $v_{e'}+\text{rad}\,Z'(\rho)$. These vectors have weights $\theta^{e'(\rho)}\rho$ and $\theta^{n-1}\rho$, respectively. This establishes (d) and (f).\qed

We define a mapping $\sigma:\R\to \R$ by
$$\sigma(\rho)=\begin{cases}\theta^{-(e(\rho)+1)}\rho& \text{if $\rho$ is nonexceptional}\\ \sigma(\rho)=\rho& \text{otherwise}.\end{cases}$$
\begin{Lem}\label{weylgroup}
\begin{itemize}
\item[(a)] $\sigma$ is a permutation of order $2m$.
\item[(b)] The orbits of $\R$ under $\sigma$ are of size $2m$ or $1$.
\item[(c)] $\sigma^{-1}(\rho)=\theta^{n-e(\rho)-1}\rho$.
\end{itemize}
\end{Lem}
\pf (a)-(b). Suppose $\rho$ is nonexceptional. By the two previous propositions $e(\sigma(\rho))+1+e(\rho)+1=n$. This equality gives\newline 
$\sigma^2(\rho)=\theta^{-(e(\sigma(\rho))+1)}\sigma(\rho)=\theta^{-n}\rho$. We conclude that $\sigma$ is a bijection and $\sigma^{2s}(\rho)=\theta^{-sn}\rho$ for every integer $s\ge 0$ hence $\sigma^{2m}=1$. Suppose that for some odd number $1+2s<2m$ $\sigma^{1+2s}(\rho)=\rho$. Then 
$$\rho=\sigma(\sigma^{2s}(\rho))=\sigma(\theta^{-sn}\rho)=\theta^{-(e(\rho)+1)}\theta^{-sn}\rho$$
because $e(\theta^{-sn}\rho)=e(\rho)$. This implies $e(\rho)+1\equiv 0\mod n$ which contradicts the fact that $e(\rho)\le n-2$.

(c) To prove (c), notice that $\sigma^{-2}(\rho)=\theta^n\rho$. So $\sigma^{-1}(\rho)=\sigma(\theta^n\rho)=\theta^{n-e(\rho)-1}\rho$ as asserted.\qed

\begin{Prop}\label{secondlink}For every nonexceptional root $\rho, \; Z'(\theta^{-e(\rho)}\rho)$ is a nonsplit extension of $L(\sigma^{-1}(\rho))$ by $L(\rho)$. When $\sigma=\sigma^{-1},\; Z'(\theta^{-e(\rho)}\rho)\ncong Z(\rho)$.
\end{Prop}
\pf Set $\zeta=\theta^{-e(\rho)}\rho$ and consider $Z'(\zeta)$. By Proposition \ref{extensions} $Z'(\zeta)$ is a nonsplit extension of $L(\theta^{n-1}\zeta)$ by $L(\theta^{e'(\zeta)}\zeta)$. Now $\zeta^{2m}=q^{-2e(\rho)}q^{e(\rho)}=q^{-e(\rho)}$ so by the definition of $e'(\zeta)$ we have $e'(\zeta)=e(\rho)$. Therefore $\theta^{e'(\zeta)}\zeta=\rho$ while $\theta^{n-1}\zeta=\sigma^{-1}(\rho)$ by part (c) of the preceding Lemma..

In the case where $\sigma=\sigma^{-1}$ (i.e. $m=1$) we have two modules with the same composition series as in the preceding result, namely $Z'(\theta^{-e(\rho)}\rho)$ and $Z'(\rho)$. These modules are nonisomorphic. For, by definition, $Z(\rho)$ has an $E$-trivial generator, while $Z'(\theta^{-e(\rho)}\rho)$ does not.\qed

We need one more result before the main theorem. This is
\begin{Lem}\label{dimensions} 
\begin{itemize}
\item[(1)] If $n$ is odd, then for every $d,\,1\le d\le n$ there are exactly $m$ nonisomorphic simple modules of dimension $d$.
\item[(2)] If $n$ is even and $\kappa=1$, then for every {\em odd} $d,\,1\le d\le n$ there are exactly $2m$ nonisomorphic simple modules of dimension $d$.
\item[(3)] If $n$ is even and $\kappa=-1$, then for every {\em even} $d,\,1\le d\le n$ there are exactly $2m$ nonisomorphic simple modules of dimension $d$.
\end{itemize}
\end{Lem}
\pf By Proposition \ref{extensions} $\dim L(\rho)=d$ if and only if $e(\rho)=d-1$. So the asserted results are a question of the number of solutions in $R_{\kappa,N}$ of the equation $\rho^{2m}=q^e$. Recall the mapping $p:R_{\kappa,N}\to R_n,\;  p(\rho)=\rho^{2m}$. For $\rho,\zeta\in\R$ $p(\rho)=p(\zeta)$ if and only if $p(\rho\zeta^{-1})=1$. Now $\rho\zeta^{-1}\in R_N$, and as $R_N$ is generated by $\theta$ with $\theta^m=q$, the image of $p$ restricted to $R_N$ is the subgroup of $R_n$ generated by $q^2$. This proves parts (1) and (2) as $p(\R)=p(R_N)$ in both cases.
% Note that if $\kappa=-1$ fixing one root, say, $\pi$ of $\R$, all other roots are of the form $\theta^r\pi$. Therefore assuming $p(\pi)=q^s$ the image $p(\R)=\{q^{2r+s}|0\le r\le n-1\}$ which equals $R_n$ because $2$ is a unit of $R_n$

Assume conditions of part (3). Pick an element $\pi\in\k$ such that $\pi^2=\theta$. Clearly $\pi\in R_{-1,N}$ and moreover every $\rho\in R_{-1,N}$ is of the form $\theta^r\pi$ for some $0\le r<N$. Now $p(\theta^r\pi)=q^{2r+1}$ which shows that $e(\rho)$ runs over all {\em odd} integers between $1$ and $n$. This proves (3).\qed

We let $P(\rho)$ denote the projective cover of $L(\rho)$. The $i$th term of the radical series of $P(\rho)$ shall be denoted by $\text{rad}^i\,P(\rho)$, and the socle of $P(\rho)$ shall be denoted by $\text{soc}\,P(\rho)$.
\begin{Thm}\label{main} 
\begin{itemize}

\item[(1)] For every nonexceptional root $\rho$
\begin{itemize}
\item[(a)] $\dim P(\rho)=2n$.
\item[(b)] The radical series of $P(\rho)$ is
$$P(\rho)\supset \text{rad}\,P(\rho)\supset\text{soc}\,P(\rho)\supset 0$$
with $\text{rad}\,P(\rho)/{\text{rad}\,^2P(\rho)}\cong L(\sigma(\rho))\oplus L(\sigma^{-1}(\rho))$.
\end{itemize}
\item[(2)] If $\rho$ is exceptional, then $P(\rho)=L(\rho)$.
\end{itemize} 
\end{Thm}
\pf Since $H$ is symmetric algebra by Theorem \ref{symalg} $\text{soc}\,P(\rho)=L(\rho)$. Let $\rho$ be a nonexceptional root. By Propositions \ref{extensions} and \ref{secondlink} there are epimorphisms of $P(\rho)$ on $Z(\rho), Z'(\theta^{-e(\rho)}\rho)$ carrying $\text{rad}\,P(\rho)$ to $L(\sigma(\rho))$ and $L(\sigma^{-1}(\rho))$, respectively. Thus $\text{rad}\,P(\rho)/{\text{rad}^2\,P(\rho)}$ has a summand $L(\sigma(\rho)\oplus L(\sigma^{-1}(\rho))$. Therefore
$$ \dim P(\rho)\ge 2\dim L(\rho)+\dim L(\sigma^{-1}(\rho))+\dim L(\sigma(\rho))=2n$$
where the last equality holds since $Z(\rho)/{L(\sigma(\rho))}=L(\rho)$\newline and $Z'(\theta^{-e(\rho)}\rho)/{L(\sigma^{-1}(\rho))}=L(\rho)$.

By general principles
\begin{equation}\label{general} H\ex=\bigoplus_{\rho\in\R}\dim P(\rho)^{\dim L(\rho)}\end{equation}
If $n$ is odd, Lemma \ref{dimensions} (1) gives
\begin{align*}\sum_{\rho\in\R}\dim P(\rho)\cdot \dim L(\rho)&=\sum_{\rho\,\text{nonexceptional}}\dim P(\rho)\cdot\dim L(\rho)\\&+\sum_{\rho\,\text{exceptional}}\dim P(\rho)\cdot \dim L(\rho)\\&\ge 2n(1+2+\ldots+ n-1)m + n^2m=Nn^2.\end{align*}
Since $\dim H\ex=Nn^2$ that expression forces $\dim P(\rho)=2n$ for all nonexceptional roots $\rho$ and $\dim P(\rho)=n^2$ for the exceptional roots. This completes the proof. The even cases are treated similarly.\qed

We denote by $B(\rho)$ the block of $H\ex$ containing $P(\rho)$. We let $\mathcal W$ denote group generated by $\sigma$.
\begin{Cor}\label{blocks} The blocks of $H\ex$ correspond to $\mathcal W$- orbits in $\R$. The correspondence is given by the formula
$$B(\rho)=\sum_{\zeta\in\mathcal{W}\rho}P(\zeta)^{\dim L(\zeta)}$$ 
\end{Cor}
\pf Recall that two projective indecomposable modules $P$ and $P'$ are {\em linked} if they share a composition factor. Linkage gives rise to an equivalence relation on the set of projective indecomposable modules. The block of $P$ is the sum of all projective indecomposable modules equivalent to $P$. By the preceding Theorem $P(\rho)$ is linked to $P(\zeta)$ if and only if $\zeta=\sigma^{\pm}(\rho)$. The formula follows.\qed
\begin{Thm}\label{quiver} Assume $\rho$ is a nonexceptional root. Then
\begin{itemize} 
\item[(a)] The Gabriel quiver $Q_{\rho}$ of the block containing $P(\rho)$ is the quiver with vertices being the isomorphism classes of simples\newline $[L(\sigma^i(\rho))],\,i\in\Bbb {Z}/{2m\Bbb Z}$ with arrows
\begin{equation*}
[L(\sigma^i(\rho))]{\raise3pt\hbox{$\displaystyle{\underset{\displaystyle\underset{b_i}\longleftarrow}{\overset{a_i}\longrightarrow}}$}} [L(\sigma^{i+1}(\rho))]
\end{equation*}
corresponding to translation by $\sigma$.
\item[(b)] The basic algebra of the block containing $P(\rho)$ is the quotient of the path algebra of $Q_{\rho}$ with relations $$a_ib_i-b_{i+1}a_{i+1}$$
and all other paths of length $\ge 2$.
\item[(c)] The block of $Q_{\rho}$ is a special biserial algebra and therefore of tame representation type.
\end{itemize}
\end{Thm}
\pf For generalities on basic algebras and their presentations by quivers with relations, the reader may consult \cite{ARS}. The present situation is similar to \cite{EGST}; however, the argument there is not compatible with the theory developed here. We proceed as follows.

As a first case assume $m>1$. By Theorem \ref{main} $\text{rad}\,P(\rho)/{\text{rad}\,^2P(\rho)}$ equals $L(\sigma(\rho))\oplus L(\sigma^{-1}(\rho))$, therefore the quiver has arrows as in the statement. Let $\k Q_{\rho}$ be the path algebra of the quiver and $\mathcal B$ be the basic algebra of $B(\rho)$. We want to construct an epimorphism of $\k Q_{\rho}$ onto $\mathcal B$ (cf. \cite [III 1.9]{ARS}) whose kernel is generated by relations in (b). Thanks to Propositions \ref{extensions} and \ref{secondlink} there are nonsplit extensions $E_i$ and $E'_i$ of $L(\sigma^i(\rho))$ by $L(\sigma^{i-1}(\rho))$  and $L(\sigma^i(\rho))$ by $L(\sigma^{i+1}(\rho))$, respectively. Since $P(\sigma^i(\rho))$ is injective both embed in $P(\sigma^i(\rho))$. Now we assign to $a_i$ an epimorphism $P(\sigma^i(\rho))\to E_{i+1}$ and to $b_i$ an epimorphism $P(\sigma^{i+1}(\rho))\to E'_i$ which we still denote by $a_i$ and $b_i$. Since $\sigma^{i-1}(\rho)\neq \sigma^i(\rho),\sigma^{i+1}(\rho)$ we have $a_i(E_i)=0$, hence $a_ia_{i-1}=0$, and likewise $b_{i-1}b_i=0$. Further $a_ib_i$ and $b_{i+1}a_{i+1}$ are endomorphisms of $P(\sigma^{i+1}(\rho))$ into its socle $L(\sigma^{i+1}(\rho))$. Therefore 
$$ \gamma_ia_ib_i=b_{i+1}a_{i+1}$$
for some $\gamma_i\in\k$ and all $i\in\Bbb {Z}/{2m\Bbb Z}$. 

We wish to change the basis to get rid of the coefficients. We will follow argumentation of \cite{EGST}. Replacing $a_i$ by $a'_i:= \gamma_i\cdots \gamma_{2m-1}a_i,\;i=1,2,\ldots,{2m-1}$, we obtain the relations
\begin{align*}a'_ib_i&=b_{i+1}a'_{i+1}\;\text{for}\;i=1,\ldots,2m-2\\
a'_{2m-1}b_{2m-1}&=b_0a_0\\
\gamma_0\cdots\gamma_{2m-1}a_0b_0&=b_1a'_1\end{align*}

We know by Theorem \ref{symalg} that $H$ is a symmetric algebra and then so are the blocks of $H\ex$ and their basic algebras. Let $\psi:\mathcal{B}\to\k$ be a symmetrizing linear form. Then 
\begin{align*}\psi(\gamma_0\cdots\gamma_{2m-1}a_0b_0)&=\psi(b_1a'_1)=\psi(a'_1b_1)\\
&=\cdots =\psi(a'_{2m-1}b_{2m-1})=\psi(b_0a_0)=\psi(a_0b_0)\end{align*}
whence $(\gamma_0\cdots\gamma_{2m-1}-1)a_0b_0\in\ker\psi$ and this element spans a one- dimensional left ideal of $H$. Thus $\gamma_0\cdots\gamma_{2m-1}=1$ as desired.
 
Moving to the case $m=1$ we let $i\in\Bbb {Z}/{2\Bbb Z}$. Now we have just two projective indecomposable modules in $B(\rho)$, call them $P_0$ and $P_1$, and let $L_i=P_i/{\text{rad}\,P_{i+1}}$. In contrast to the previous case, there doesn't exist a uniform choice of length two indecomposable submodules in $\text{rad}\,P_i$. We argue as follows. As above by Proposition \ref{secondlink} there are noisomorphic submodules $E_1$ and $E_1^{\prime}$
in $\text{rad}\,P_1$ with the top composition factor $L_0$. Consequently, there are epimorphisms $a_0:P_0\to E_1$ and $b_1:P_0\to E_1^{\prime}$. Put $F_0=\ker a_0$ and $F_0^{\prime}=\ker b_1$. Both $F_0$ and $F_0^{\prime}$ lie in $\text{rad}\,P_0$ hence by Theorem \ref{main} they are extensions of $L_0$ by $L_1$. Therefore there are epimorphisms $a_1:P_1\to F_0$ and $b_0:P_1\to F_0^{\prime}$; hence $a_0a_1=0=b_1b_0$ holds by construction. However, since $a_i,b_i$ are nonisomorphisms they send the socle of $P_j$ to zero for all $i$. Therefore every path of length $\ge 2$ is either zero or spans a $1$-dimensional ideal of $\mathcal B$. Since $\mathcal B$ is  symmetric $a_1a_0=0=b_0b_1$ as well, whence $E_1\subset\ker a_1$ and $E_1^{\prime}\subset\ker b_0$. As $\dim\ker a_1=n=\dim E_1$ we see that $E_1=\ker a_1$ and similarly $E_1^{\prime}=\ker b_0$. From this we conclude that $F_0\neq F_0^{\prime}$. Indeed, since $\dim\text{Hom}(P_1,F_0)=[F_0:L_1]=1$ were $F_0=F_0^{\prime}$ we would have $a_1=\alpha b_0,\alpha\in\Bbbk^{\bullet}$, which forces $E_1=E_1^{\prime}$, a contradiction. It follows that $a_0b_0$ and $b_1a_1$ are nonzero homomorphisms $P_1\to\text{soc}\,P_1$, and similar conclusions hold for $b_0a_0$ and $a_1b_1$. Therefore there are $\gamma_i\in\Bbbk^{\bullet}$ such that 
\begin{align*}\gamma_1a_0b_0&=b_1a_1\\\gamma_2b_0a_0&=a_1b_1\end{align*}
But then the calculation
\begin{align*}\psi(b_1a_1)&=\gamma_1\psi(a_0b_0)=\psi(a_1b_1)\\&=\gamma_2\psi(b_0a_0)=\gamma_2\psi(a_0b_0)\end{align*}
gives $\gamma_1=\gamma_2$. Replacing $a_0$ by $\gamma_1a_0$ we obtain the desired result.\qed

The proof of (c) follows directly from the definition and theory of special biserial algebras, see \cite{Erd}.\qed

Lastly, we want to describe embeddings of blocks and projective indecomposable modules in $H$. First, we take the subalgebra $A_{\lambda}$ of $H\ex$ of section \ref{casimir}. From its definition we have $A_{\lambda}\cong \k[t]/{(t^N-\kappa)}$, hence $A_{\lambda}=\oplus \k e_{\zeta}$ where $e_{\zeta}$ is a primitive idempotent uniquely determined by the property $Ke_{\zeta}=\zeta e_{\zeta}$. Secondly, we use primitive idempotents $\epsilon_{\mu}$ associated to roots $\mu\in R_{\kappa,n}$ defined in section \ref{babyverma}. For every $\rho\in R_{\kappa,N},\rho^m\in R_{\kappa,n}$ hence gives rise to the idempotent $\epsilon_{\rho^m}$. Abusing notation we write the latter as $\epsilon_{\rho}$.
\begin{Prop}\label{blockidempotents} Suppose $\rho$ is a nonexceptional root. Then
\begin{itemize}
\item[(1)] $B(\rho)=H\ex\epsilon_{\rho}$.
\item[(2)] Every nonsimple projective indecomposable module of $H$ is isomorphic to $H\ex\epsilon_{\rho}e_\zeta$ for some $\lambda,\rho,\zeta$.
\end{itemize}
\end{Prop}
\pf It is elementary to see that for a root $\mu$ of the minimal polynomial $f_{\lambda}(t)$, $\epsilon_{\rho}(\mu)=1$ if $\mu=-\eta'D(\rho)$, and zero, otherwise. For every simple $H\ex$- module $L(\zeta)$ $C_{\lambda}$ acts on $L(\zeta)$ via multiplication by $-\eta'D(\zeta)$. Therefore $\epsilon_{\rho}$ acts on $L(\zeta)$ via multiplication by $\epsilon_{\rho}(-\eta'D(\zeta))$. Now $D(\zeta)=D(\rho)$ if and only if $\zeta^m=\rho^m$ or $\zeta^m\rho^m=q^{-1}$. Since $\rho$ is nonexceptional, $q^{-1}\rho^{-m}\neq \rho^m$; hence both cases occur. Next observe that for every fixed element $\mu\in R_{\kappa,n}$ the set of solutions to the equation $x^m=\mu$ in $R_{\kappa,N}$ has $m$ elements. Therefore $Y=\{\zeta| D(\zeta)=D(\rho)\}$ has $2m$ elements.

On the other hand by Lemma \ref{weylgroup}, $\zeta\in\mathcal{W}\cdot\rho$ if and only if $\zeta=\theta^{-np}\rho$ or $\zeta=\theta^{-np-e(\rho)-1}\rho$. In the first case $\zeta^m=\rho^m$ and in the second $\zeta^m\rho^m=q^{-1}$. Thus $\mathcal{W}\cdot\rho\subset Y$, hence by Lemma \ref{weylgroup} (b) they are equal. It follows that $\epsilon_{\rho}$ acts as the identity on $B(\rho)$ and annihilates every other block, whence (1).

The number of projective indecomposable modules in a decomposition of $B(\rho)$ equals $\sum_{\zeta\in\mathcal{W}\cdot\rho}\dim L(\zeta)$. For every $\zeta\in\mathcal{W}\cdot\rho$, $\sigma(\zeta)\neq\zeta$, hence there are $m$ distinct pairs $\{\zeta,\sigma(\zeta))$ in $\mathcal{W}\cdot\rho$. For each such pair Proposition \ref{extensions} (c) gives $\dim L(\zeta)+\dim L(\sigma(\zeta))=n$. Thus $B(\rho)$ is the direct sum $N$ indecomposable summands. Since we have $N$ idempotents $e_{\zeta}$, the result follows.\qed

We turn to the case of a projective simple module $P(\rho)$. 
\begin{Prop}\label{simpleproj} Suppose $\rho$ is an exceptional root. Let $f=e_{\rho}E^{n-1}F^{n-1}$. Then 
$$H\ex f\cong P(\rho)$$
and $H\ex f\;\text{ is a direct summand of}\quad H\ex$.
\end{Prop}
\pf Evidently $f$ has weight $\rho$. Further, it is elementary to derive the formula $e_{\rho}=\displaystyle\frac{\rho}{N\kappa}\prod_{\zeta\neq\rho}(K-\zeta)$ and then a simple calculation gives $Ee_{\rho}=e_{\theta\rho}E$, hence $f$ is $E$-trivial, which proves the first assertion.

Next we adjust Kac's formula \cite[(1.3.1)]{DeCK} to our context obtaining
\begin{equation}\label{kac}[E^s,F^r]=\eta(\sum_{i=1}^{\text{min}(r,s)}F^{r-i}H_i^{s,r}E^{s-i})\end{equation}
where $H_i^{s,r}=(r)_q\cdots (r-i+1)_q{\binom si}_q\prod_{j=1}^{j=i}(K^{-m}-q^{i+j-r-s}K^m)$. We now compute $f^2$ using \eqref{kac}. Since $E^n=0$ we deduce
$$f^2=e_{\rho}H_{n-1}^{n-1,n-1}f=cf$$
where $c=\eta(n-1)_q!\prod_{j=1}^{n-1}(\rho^{-m}-q^{j-n+1}\rho^m)$ and since $\rho$ is exceptional, every factor in $c$ is nonzero, hence $c\neq 0$. Thus $c^{-1}f$ is an idempotent, and the proof is complete.\qed

\end{document}